\documentclass[12pt,reqno]{amsart}
\usepackage{amsthm,amsfonts,amssymb,amscd}
\usepackage[all]{xy}

\newtheorem{theorem}[subsection]{Theorem}

\newtheorem{proposition}[subsection]{Proposition}

\newtheorem{lemma}[subsection]{Lemma}
\newtheorem{corollary}[subsection]{Corollary}

\theoremstyle{definition}

\newtheorem{proposition-definition}[subsection]{Proposition-Definition}

\theoremstyle{remark}
\newtheorem{remark}[subsection]{Remark}

\newcommand{\fatdot}{{\scriptscriptstyle \bullet}}

\newcommand{\dual}{{{\scriptscriptstyle \vee}}}

\newcommand{\Ext}{\operatorname{Ext}\nolimits}

\newcommand{\gr}{\operatorname{gr}}
\newcommand{\Gr}{\operatorname{Gr}}
\newcommand\Hilb{{\operatorname{Hilb}\nolimits}}

\newcommand{\Hom}{\operatorname{Hom}\nolimits}
\newcommand{\id}{\operatorname{id}\nolimits}
\newcommand{\im}{\operatorname{im}\nolimits}

\newcommand{\Pic}{\operatorname{Pic}\nolimits}

\newcommand{\Res}{\operatorname{Res}}
\newcommand{\rk}{\operatorname{rk}\nolimits}

\newcommand\Sym{{\operatorname{Sym}\nolimits}}
\newcommand{\Tr}{\operatorname{Tr}}

\newcommand{\CC}{{\mathbb C}}

\newcommand{\ZZ}{{\mathbb Z}}
\newcommand{\QQ}{{\mathbb Q}}

\newcommand{\HHB}{{\mathbf H}}

\newcommand{\OOO}{{\mathcal O}}

\newcommand{\JJJ}{{\mathcal J}}

\newcommand{\EEE}{{\mathcal E}}

\newcommand{\HHH}{{\mathcal H}}

\newcommand{\FFF}{{\mathcal F}}

\newcommand{\DDD}{{\mathcal D}}
\newcommand{\KKK}{{\mathcal K}}

\newcommand{\MMM}{{\mathcal M}}

\newcommand{\TTT}{{\mathcal T}}
    
\newcommand{\XXX}{{\mathcal X}}

\newcommand{\VVV}{{\mathcal V}}
\newcommand{\WWW}{{\mathcal W}}
\newcommand{\YYY}{{\mathcal Y}}

\newcommand\alp{\alpha}

\newcommand\eps{\epsilon}

\newcommand{\into}{\hookrightarrow}

\newlength{\rrrr}
\newcommand{\isom}[1]{{\settowidth{\rrrr}{$\scriptstyle{x#1x}$}
\xrightarrow{\makebox[\rrrr]{$\scriptstyle{#1}$}}
\hspace{-0.5\rrrr }\hspace{-1.1 em}
\raisebox{- 0.5 ex}{$\sim$}\hspace{0.7\rrrr }
}}
\newcommand{\isoto}{{\lra\hspace{-1.3 em}
\raisebox{ 0.6 ex}{$\textstyle\sim$}\hspace{0.8 em}}}
\newcommand\lra{{\longrightarrow}}

\newcommand\rar{\rightarrow}

\renewcommand{\bar}[1]{\overline{#1}}

\author{D. Markushevich}

\address{D. M.: Math\'ematiques - b\^{a}t. M2, Universit\'e Lille 1,
F-59655 Villeneuve d'Ascq Cedex, France}
\email{markushe@math.univ-lille1.fr}


\title{An integrable system of K3-Fano flags}

\begin{document}
\begin{abstract}
Given a K3 surface $S$, we show that the relative intermediate
Jacobian of the universal family of
Fano 3-folds $V$ containing $S$ as an anticanonical divisor
is a Lagrangian fibration.
\end{abstract}
\maketitle

\section*{Introduction}

Beauville \cite{B2} defined 
the moduli stack $\FFF^R_g$ of K3-Fano flags $S\subset V$, where $S$ is
a K3 surface and $V$ an $R$-polarized Fano 3-fold 
containing $S$ as an anticanonical divisor.
Here $R$ is an even integral lattice of signature $(1, \rk (R)-1)$
with a distinguished element $\rho$, and an $R$-polarization of a Fano variety
$V$ is an isomorphism $R\rar \Pic(V)$ which sends $\rho$ to $-K_V$ and
such that the bilinear form on $R$ is the pullback of the
form $(D_1,D_2)\mapsto (-K_V\cdot D_1\cdot D_2)$ on $\Pic(V)$.
The subscript $g$ is an integer parameter, defined by $2g-2=(-K_{V})^3$.
Beauville also showed that the natural map
$s^R_g:\FFF^R_g\rar\KKK^R_g$ to the moduli stack $\KKK^R_g$ of $R$-polarized K3 surfaces with a very ample class of degree $2g-2$ is generically surjective.

The relative dimension of $s^R_g$ at $(S,V)$ is $b_3(V)/2$,
and the fiber $\FFF^R_S:=(s^R_g)^{-1}(S)$ is a smooth stack.
Let $\VVV\rar\FFF^R_S$ be the universal family of Fano 3-folds
over $\FFF^R_S$ and $\JJJ=J(\VVV/\FFF^R_S)$ its relative intermediate
Jacobian. In this note, we prove that $\JJJ$ is holomorphically
symplectic and that the structure morphism $\JJJ\rar\FFF^R_S$ is
a Lagrangian fibration.

Our proof is Hodge theoretic. Some particular cases have been treated earlier
by different geometric methods.
Hassett--Tschinkel (\cite{HT}, Proposition 7.1) defined a Lagrangian fibration on the length-2 punctual Hilbert scheme $S^{[2]}$ of a K3 surface $S$
which is a complete intersection of three quadrics. The Fano threefolds
$V$ containing $S$ are intersections of 2 quadrics, 
and the fibers of the Lagrangian
fibration are the Fano surfaces $F(V)$ parametrizing lines in $V$. It follows
from \cite{N}, \cite{NR} that $F(V)\simeq J(V)$, the intermediate Jacobian of $V$.
Thus the Lagrangian fibration of Hassett--Tschinkel is a torsor under
the relative intermediate Jacobian $\JJJ$, which implies that
$\JJJ$ is a Lagrangian fibration itself.

Beauville constructed  in \cite{B1} a Lagrangian fibration
associated to the complete family of cubic 3-folds $V_t$ containing
a given K3 surface $S$ which is a complete intersection of a quadric and a cubic; here $g=13$ and $R=\ZZ\cdot \frac{1}{2}\rho$. The phase space $\MMM$ is a family of moduli spaces $M_t$ of vector bundles
on the cubic 3-folds $V_t$ containing $S$, and each $M_t$ is isomorphic to
an open subset of the intermediate Jacobian $J(V_t)$. Thus it is
an open part of a torsor under the relative intermediate Jacobian
$\JJJ\rar\FFF^R_S$.

The authors of \cite{IM} introduced a Lagrangian
structure on $J(\VVV/\FFF^R_S)$ for $3\leq g\leq 10$ over an open subset
of $\FFF^R_S$, where $S$ is assumed generic and $R= \ZZ\rho$. Their construction
is a generalization of that of \cite{HT}: $J(\VVV/\FFF^R_S)$ is identified
with the relative Picard variety of a family of Lagrangian
subvarieties of $S^{[2]}$, and the relative Picard variety is itself
Lagrangian by a result of Donagi--Markman \cite{DM}.

The approaches of \cite{B1} and \cite{IM} are based upon
the observation of Tyurin that for a K3-Fano flag
$(S,V)$, the map sending stable sheaves on $V$ to their restrictions to $S$
under certain hypotheses embeds the moduli spaces of sheaves on $V$ as 
Lagrangian subvarieties of moduli spaces on $S$. Tyurin \cite{Tyu}
originally stated it for moduli of stable vector bundles,
and Thomas \cite{Th} extended to the restriction map
$\Hilb^{t\mapsto dt+c} (V)\rar S^{[d]}$
defined on sufficiently nice degree-$d$ curves in $V$;
the authors of \cite{IM} reproved the result of Thomas and
applied it to conics, that is for $d=2$.

The referee of the present paper evoked yet another way to construct a
Lagrangian fibration from K3-Fano flags: take a pair of K3-Fano
flags $(S,V_1)$, $(S,V_2)$ with the same $S$, and consider the normal
crossing variety $X=V_1\cup V_2$, where $V_1$, $V_2$ intersect transversely
along $S$. Under certain hypotheses, $X$
can be smoothed to a Calabi--Yau threefold. Then the
Donagi--Markman integrable system over the moduli space of Calabi--Yau
threefolds obtained by smoothing $X$ defines the ``boundary''
integrable system whose Liouville tori are products $J(V_1)\times J(V_2)$.
A similar construction was described in \cite{Do}, Section 6, where
the normal crossing Calabi--Yau was the union of two threefolds,
each an elliptic fibration over a del Pezzo surface rather than a Fano
threefold. Our Lagrangian fibration seems to be a direct factor of the
boundary integrable system.

In Sect. \ref{DMC}, we recall the Donagi--Markman criteria for the
existence of a quasi-Lagrangian
and Lagrangian structures. In Sect. \ref{VHS},
we gather basic facts on the variations of Hodge structures necessary
to verify the weak cubic condition of Donagi--Markman. We verify it
in Sect. \ref{K3Ff}, which provides $\JJJ$ with a quasi-Lagrangian
structure. To get a genuine Lagrangian structure, we need variations
of mixed Hodge structures, which are sketched in Sect. \ref{VMHS}.
Finally, in Sect. \ref{SS}, we verify the sufficient condition of Donagi--Markman
for the existence of a Lagrangian structure on $\JJJ$.

{\sc Acknowledgements.} The author acknowledges with pleasure the hospitality of
the Max-Planck-Institut f\"ur Mathematik in Bonn, where was done
the work on the present paper. He thanks the referee for interesting comments
and for his remark suggesting that the local Lagrangian structure defined
in the first version of the paper should be in fact global; this remark prompted
the author to improve the main result.

\section{Donagi--Markman cubic}\label{DMC}

We will start by a reminder on the Kodaira-Spencer map. Let $f:\VVV\rar B$
be a family of compact complex manifolds. The space of infinitesimal deformations
of any special fiber $V=V_b$ of $f$ ($b\in B$) over the Artinian
ring $\CC [\eps]/(\eps^2)$
is identified
with $H^1(V,\TTT_V)$. Any tangent vector $a\in T_bB$ defines
such an infinitesimal deformation, so we have the natural
map $KS(f)_b:T_bB \rar H^1(V,\TTT_V)$, called Kodaira--Spencer map.
Under the assumption that $h^1(\TTT_{V_b})$ is constant on $B$, the maps $KS(f)_b$
fit into the bundle map $KS(f):\TTT_B\rar R^1f_*\TTT_{\VVV/B}$,
where $\TTT_{\VVV/B}$ denotes the vertical tangent bundle of $f$.
One can also interpret $KS(f)$ as the image of
the extension class of the tangent bundle sequence
\begin{equation}\label{tbs1}
0\lra \TTT_{\VVV/B}\lra \TTT_{\VVV}\lra f^*\TTT_B\lra 0
\end{equation}
in $\Ext^1(f^*\TTT_B, \TTT_{\VVV/B})=
H^1(\VVV,f^*\TTT_B^\dual\otimes\TTT_{\VVV/B})$ under the natural map
$$
H^1(\VVV,f^*\TTT_B^\dual\otimes\TTT_{\VVV/B})\lra
H^0(B,\TTT_B^\dual\otimes R^1f_*\TTT_{\VVV/B}).
$$

Replacing the $V_b$ by their intermediate Jacobians $J^k(V_b)$,
we obtain a family $\pi=J^kf:\XXX\rar B$ of complex tori,
and its Kodaira--Spencer class $KS(\pi)$ is related to $KS(f)$
by means of the Gauss--Manin connection, see Sect. \ref{VHS}.

Given a family of abelian varieties or complex tori $\pi :\XXX\rar B$,
we will say that it is quasi-Lagrangian if $\dim B=\frac{1}{2}\dim\XXX$ and
there exists a nondegenerate holomorphic
2-form $\omega$ on $\XXX$
such that the fibers of $\pi$ are maximal isotropic with respect to $\omega$.
If such an $\omega$ exists only locally over $B$, we say that $\pi$ 
is locally quasi-Lagrangian, or that it has a local quasi-Lagrangian structure.
Requiring in the above definitions that $\omega$ is closed, and hence symplectic,
we obtain the notions of a Lagrangian and local Lagrangian structures.
In \cite{DM} Donagi and Markman determine conditions on $\pi$
which are necessary and/or sufficient for the existence a 
\mbox{(quasi-)}Lagrangian structure. 
Let $\EEE$ denote the direct image $\pi_* \TTT_{\XXX/B}$ on $B$ of the vertical tangent bundle
of $\XXX$. As the fibers of $\pi$ are isotropic, $\omega$
induces an isomorphism
$
j:\TTT_B^\dual\rightarrow \EEE.
$
As before, the tangent bundle sequence
\begin{equation}\label{tbs2}
0\lra \TTT_{\XXX/B}\lra \TTT_\XXX\lra\pi^*\TTT_B\lra 0
\end{equation}
determines the class
$$
KS(\pi)\in H^0(B,\TTT_B^\dual\otimes R^1\pi_*\TTT_{\XXX/B})=
H^0(B,\TTT_B^\dual\otimes R^1\pi_*\OOO_\XXX\otimes \EEE).
$$
We will say that a complex torus $\CC^g/\Lambda$ is quasi-polarized
if $\Lambda$ has a skew-symmetric bilinear form satisfying
the first Riemann bilinear relation (the symmetry of the period
matrix), but not necessarily the second one (positive definiteness
of the imaginary part of the period matrix).
Assume that $\XXX/B$ is a family of quasi-polarized complex tori.
Then the quasi-polarization defines an isomorphism
$\theta:R^1\pi_*\OOO_\XXX\rar\EEE$. 
The Donagi--Markman cubic associated to $(\pi,\theta,j)$
is the element
$Q\in H^0(B,\EEE^{\otimes 3})$ defined by $Q=H^0(j\otimes\theta\otimes \id_\EEE)(KS(\pi))$.

\begin{lemma}[Donagi--Markman, \cite{DM}]\label{weak}
Let $\pi:\XXX\rar B$ be a family of quasi-polarized $g$-dimensional
complex tori over a complex manifold $B$ of dimension $g$ and
$j:\TTT_B^\dual\rightarrow \EEE$ an isomorphism.
Then there exists
a local quasi-Lagrangian structure $\omega$ on $\XXX/B$
inducing the given isomorphism
$j:\TTT_B^\dual\rightarrow \EEE$ if and only if
the associated Donagi--Markman cubic is symmetric,
that is $Q\in H^0(B,\Sym^3(\EEE))$. If, moreover, $\pi$ admits
a cross-section $O$, then the symmetry of $Q$ guarantees the existence
of a global quasi-Lagrangian structure $\omega$ on $\XXX/B$ inducing $j$, and
$\omega$ is unique if we impose the additional condition that $O$
is isotropic with respect to $\omega$, that is $\omega|_O\equiv 0$.
\end{lemma}

This lemma gives a necessary condition for the
existence of a Lagrangian structure on $\XXX/B$.
Now we will state a sufficient condition.
Let $\HHH_k(\XXX/B,\ZZ)$ denote the local system of $k$-th homology
groups of fibers of $\pi$ with integer coefficients. The integration of \mbox{1-forms} over \mbox{1-cycles} on the
fibers of $\pi$ defines a natural embedding $\HHH_1(\XXX/B,\ZZ)
\into \mbox{$\pi_*(\Omega^1_{\XXX/B})^\dual\simeq \EEE$}$.

\begin{lemma}[Donagi--Markman, \cite{DM}]\label{strong}
Under the assumptions of Lemma \ref{weak}, 
assume that for any open subset $U$ of $B$ and for any relative
$1$-cycle $\gamma\in H^0(U,\HHH_1(\XXX/B,\ZZ))$,
the preimage $j^{-1}(\gamma)$ in $H^0(U,\TTT_B^\dual)$ is a
closed $1$-form on $U$. Then $\pi$ possesses a local Lagrangian structure
inducing $j$.
If, moreover, $\pi$ admits
a cross-section $O$, then there exists a unique global Lagrangian structure $\omega$ on $\XXX/B$, inducing $j$ and such that $O$ is Lagrangian
with respect to $\omega$.
\end{lemma}

\section{Variations of Hodge structures}\label{VHS}

Let $f:\VVV\rar B$ be a family of smooth projective or
compact K\"ahler manifolds
of dimension~$n$. Then the vector bundles $\HHH^k=R^kf_*\CC\otimes_\CC\OOO_B$
carry the structure of a variation of Hodge structures (VHS) of weight $k$. 
This means the following:
\begin{enumerate}
\item $\HHH^k$ contains a local system $\HHH^k_\ZZ$ of integral
lattices such that $\HHH^k=\HHH^k_\ZZ\otimes_\ZZ\OOO_B$. This, in particular, implies the
existence of a natural real structure, defining the complex
conjugation $s\rar \bar s$ on the sections of $\HHH^k$, and
the existence of a flat connection $\nabla:\HHH^k\rar \HHH^k\otimes\Omega^1_B$,
called Gauss--Manin connection. The latter is defined by the requirement
that the sections of $\HHH^k_\ZZ$ are flat.
\item $\HHH^k$ possesses a decreasing filtration
$\HHH^k=F^0\HHH^k\supset F^1\HHH^k\supset \ldots \supset F^k\HHH^k\supset
F^{k+1}\HHH^k=0$ by holomorphic subbundles such that
$F^i\HHH^k\cap \bar{F^{k-i-1}\HHH^k}=0$.
\item The Griffiths transversality condition:
$\nabla (F^i\HHH^k)\subset F^{i-1}\HHH^k\otimes\Omega^1_B$.
\end{enumerate}

For each fiber $V_b=V$ of $f$, 
$$\HHH^k|_b=H^k(V,\CC),\ \ 
F^i\HHH^k|_b=\bigoplus\limits_{p\geq i}H^{p,k-p}(V,\CC),\ \ 
\HHH^k_\ZZ|_b=H^k(V,\ZZ)/\mbox{(torsion)}.
$$

The Gauss-Manin connection is only $\CC$-linear, but when restricted to $F^k$, it
induces an \mbox{$\OOO_B$-linear} map $\bar\nabla:F^k\rar (F^0/F^k)\otimes \Omega^1_B$,
whose image is in fact contained in $(F^{k-1}/F^k)\otimes \Omega^1_B$.
By \cite{Gr}, the graded version of the Gauss-Manin connection
$\gr^i\nabla:F^i/F^{i+1}\rar F^{i-1}/F^{i}\otimes\Omega^1_B$
acts by the contraction with
the Kodaira--Spencer class $KS(f)\in H^0(B,\TTT_B^\dual\otimes 
R^1f_*\TTT_{\VVV/B})$. Explicitly,
we can use the Dolbeault isomorphism to identify $F^i\HHH^k|_b$
with $\bigoplus_{p\geq i}H^{k-p}(V,\Omega^p)$, and for
each $a\in T_bB$, the graded covariant derivative $\gr^p\nabla_a$
in the direction of $a$
acts on $\xi\in H^{q}(V,\Omega^p)$ by the formula
\begin{equation}\label{nabla}
\gr^p\nabla_a:
\xi\mapsto (KS(f)|_b(a))\lrcorner\, \xi\in H^{q+1}(V,\Omega^{p-1}) , \ \ 
H^1(V,\TTT_V)\times H^{q}(V,\Omega^p)\xrightarrow{\ \lrcorner\ }
H^{q+1}(V,\Omega^{p-1}).
\end{equation}

To a VHS 
$(\HHH^{2k-1},\HHH^{2k-1}_\ZZ,F^\fatdot)$ of odd weight ${2k-1}$,
one can associate a
family of complex tori $\pi:\JJJ^k\rar B$
of relative dimension $\frac{1}{2}\rk\HHH^{2k-1}$:
$$
\JJJ^k=\HHH^{2k-1}/(\HHH^{2k-1}_\ZZ+F^k\HHH^{2k-1})
$$
If $\HHH^{2k-1}$ comes from a morphism $f$ as above, $\JJJ^k=
\JJJ^k(\VVV/B)$ is called the relative $k$-th intermediate Jacobian
of $f$. In this case we can also write
$$
\JJJ^k\cong (F^{n-k+1}\HHH^{2n-2k+1})^\dual/\HHH_{2n-2k+1}(\VVV/B,\ZZ).
$$

If $n=2k-1$, then the $k$-th intermediate Jacobian
has a natural quasi-polarization, given by the intersection form on $\HHH^{2k-1}_\ZZ$. If in addition the Hodge structure 
on $\HHH^{2k-1}$ is of height 1, that is
$F^{k+1}=0$, $F^{k-1}=\HHH^{2k-1}$, then
it also satisfies the second Riemann bilinear relation
and hence is a polarization.

From now on, we are assuming that $n=2k-1$. Let $\XXX=\JJJ^k$. We are going
to describe the extension class $KS(\pi)$ of (\ref{tbs2}) in terms of the
Gauss--Manin connection $\bar\nabla:F^k\rar (F^0/F^k)\otimes \Omega^1_B$.
The natural pairing $F^k\times (F^0/F^k)\rar \HHH^{2n}\simeq\CC$ identifies
$F^0/F^k$ with the dual of $F^k$, and
$F^0/F^k\cong \pi_*\TTT_{\XXX/B}$. As in the previous
section, we
will denote this bundle by $\EEE$. Thus we can interpret $\bar\nabla$
as a map $\bar\nabla:\TTT_B\rar (\EEE)^{\otimes 2}$. 

\begin{lemma}\label{bar-nabla}
The extension class of the tangent bundle sequence of the map
$\pi:\XXX=\JJJ^k\rar B$ coincides with $\bar\nabla$.
\end{lemma}

\begin{proof}
The VHS of weight 1 of $\pi$ is easily reconstructed from the VHS of
weight $2k-1$ of $f$. The underlying vector bundle
$\HHH^1(\XXX/B,\CC)$ is just $\HHH^{2k-1}(\VVV/B,\CC)$,
similarly for integer lattices we have $\HHH^1(\XXX/B,\ZZ)=\HHH^{2k-1}(\VVV/B,\ZZ)$,
and the Hodge filtration is given by 
$F^1\HHH^1(\XXX/B,\CC)=F^k\HHH^{2k-1}(\VVV/B,\CC)$, $F^2\HHH^1(\XXX/B,\CC)=0$.
Hence
the Gauss--Manin connection on $\HHH^1(\XXX/B,\CC)=\HHH^{2k-1}(\VVV/B,\CC)$ constructed from $\pi$ coincides
with that constructed from $f$. By the cited result of Griffiths, applied
to $\pi$, the Gauss--Manin connection is given by the contraction
with $KS(\pi)$. The component 
$\bar\nabla:F^1\HHH^1\simeq \EEE^\dual \rar \HHH^1/F^1\HHH^1\simeq\EEE\otimes\Omega^1_B$ determines
completely $KS(\pi)$, and coincides with it if interpreted as a
map $\TTT_B\rar\EEE^{\otimes 2}$.
\end{proof}

We will apply the above lemma in the case when
the fibers of $f$ are Fano 3-folds. Then $H^{3,0}(V_b)=0$ for all $b\in B$,
so the Hodge structure on $\HHH^3$ is of height 1 and $\JJJ=\JJJ^1(\VVV/B)$
is a family of polarized abelian varieties of dimension $b_3(V_b)/2$.

\section{Deformations of K3-Fano flags}\label{K3Ff}

Recall that a Fano 3-fold is by definition a 3-dimensional nonsingular
projective variety $V$ with ample anticanonical divisor $-K_V$.
Iskovskih \cite{I1}, \cite{I2} classified Fano 3-folds with Picard
number 1, and Mori--Mukai classified all the remaining ones \cite{MM};
see also \cite{IP}.
There are 105 deformation classes of Fano 3-folds,
and around half of them have $h^{2,1}=0$, so that there is no
integrable system of intermediate Jacobians associated to them.
For the sequel, fix some class of Fano 3-folds $V$ with $h^{2,1}\neq 0$ and
consider a K3-Fano flag $(S,V)$. This means that $S$ is a K3 surface
embedded in $V$ as an anticanonical divisor. Let $R$ be the Picard lattice of $V$,
$(-K_V)^3=2g-2$; by Lefschetz Theorem, $R$ is embedded into $\Pic(S)$.
Denote by $i$ the natural embedding $S\into V$.

According to \cite{Ka}, the deformation theory of a pair $D\subset X$,
consisting of a connected compact manifold $X$ and a normal crossing
divisor $D$ in it, is governed by the sheaf
$$\TTT_X(-\log D)=\{ v\in \TTT_X\mid vI_D\subset I_D\},$$
where $I_D$ denotes the ideal sheaf of $D$.
Beauville \cite{B2} gives a proof of this fact in the algebraic setting. 
Applying it to the moduli stack
of K3-Fano flags $\FFF^R_g$, we obtain that
the tangent space to $\FFF^R_g$ at $(S,V)$
is canonically isomorphic to $H^1(V,\TTT_V(-\log S))$, and the obstruction
space $H^2(V,\TTT_V(-\log S))$ is zero, so the first order
deformations of $(S,V)$ are unobstructed. 

The space of infinitesimal
deformations of the complex structure on $S$ is $H^1(S,\TTT_S)$,
and as $H^2(S,\TTT_S)=0$, the infinitesimal
deformations are also unobstructed.
Let $\sigma:R\rar H^1(S,\Omega^1_S)$ be the natural map, defined
as the composition
$R\isoto\Pic(V)\xrightarrow{i^*}\Pic(S)\xrightarrow{c_1}H^1(S,\Omega^1_S)$,
where $c_1$ denotes the first Chern class. Fixing some generator $\alp_S
\in H^0(S,\Omega^2_S)\simeq \CC$, we have the contraction isomorphism
$\TTT_S\isom{\cdot\lrcorner\alp_S}\Omega^1_S$, which we can use
to identify $H^1(S,\TTT_S)$ with $H^1(S,\Omega^1_S)$. Under
this identification, the tangent space $T_{[S]}\KKK^R_g$
is the orthogonal complement to $\sigma(R)$ with respect
to the natural bilinear form on $H^1(S,\Omega^1_S)$:
$T_{[S]}\KKK^R_g=\sigma(R)^\perp\subset H^1(S,\TTT_S)$.

Further, Beauville shows that the image of the natural map
$H^1(r):H^1(V,\TTT_V(-\log S))\rar H^1(S,\TTT_S)$
induced by the restriction $r:\TTT_V(-\log S) \rar \TTT_S$
coincides with $\sigma(R)^\perp$ and that $H^1(r)$ is
the differential of the forgetful morphism
$s^R_g:\FFF^R_g\rar \KKK^R_g$ at $[(S,V)]$. This implies, in particular,
that $s^R_g$ is a submersion. Hence the fiber $\FFF^R_S:= (s^R_g)^{-1}([S])$
is a smooth stack whose tangent space at $[(S,V)]$ is  $\ker H^1(r)$.

To identify $\ker H^1(r)$, we follow Beauville who considers
the natural exact triple
\begin{equation}\label{log-tbs}
0\lra \TTT_V(-S)\lra\TTT_V(-\log S)\xrightarrow{\ r\ }\TTT_S \lra 0.
\end{equation}
Choose a generator $\alp_V\in H^0(V,\Omega^3_V(\log S))\simeq \CC$.
Then we have the contraction isomorphism
$\TTT_V(-S)\isom{\cdot\lrcorner\, \alp_V}\Omega^2_V$. The
long exact cohomology sequence associated to (\ref{log-tbs}) acquires the form
$$
0\rar H^1(V,\Omega^2_V)\rar
H^1(V,\TTT_V(-\log S))\xrightarrow{H^1(r)}
H^1(S,\TTT_S)\rar
H^2(V,\Omega^2_V)\rar 0
$$
and shows that $T_{[(S,V)]} \FFF^R_S=\ker(H^1(r))=H^1(V,\Omega^2_V)$.
In particular, $\dim \FFF^R_S =b_3(V)/2$.

\begin{proposition}\label{ebQ}
Let $R$, $S$, $\FFF^R_S$ be as above. Let $f=\VVV\rar B=\FFF^R_S$ be
the universal family of Fano $3$-folds having $S$ as an anticanonical
divisor and $\pi:\JJJ=\JJJ(\VVV/B)\rar B$ its relative intermediate
Jacobian. Let $b=[(V,S)]$ be a point of $B$. 
Fix a generator $\alp_V$ of $H^0(V,\Omega^3_V(\log S))$ and denote by 
$\alp$ the corresponding isomorphism $\cdot\lrcorner\alp_V:
\TTT_V(-S)\isoto \Omega^2_V$. Let $h_b=H^1(\alp)$ denote 
the isomorphism $T_bB=H^1(\TTT_V(-S))\isoto H^1(V,\Omega^2_V)$ induced by
$\alp$. 

The following statements hold:

(i) The value $KS(f)_b$ of the Kodaira--Spencer class of $f$ at $b$ is the
natural map $H^1(s):T_bB=H^1(\TTT_V(-S))\rar H^1(\TTT_V)$
associated to the morphism $s$ in the
exact triple
$$
0\lra T_V(-S)\xrightarrow{\ s\ }T_V\lra T_V|_S\lra 0.
$$

(ii) The tangent space at $0$ to the
intermediate Jacobian $J^2(V)$ being  $H^1(V,\Omega^2_V)^\dual
\cong H^2(V,\Omega^1)$, we define the Donagi--Markman cubic $Q$ of $\pi$ at $b$
with the help of the isomorphism $j_b=(^th_b)^{-1}:T_bB^\dual\rar T_0J^2(V)$.
Then $Q$, as a cubic form on $H^1(V,\Omega^2_V)$, coincides
with the composition of the cup-product
$H^1(V,\Omega^2_V)^{\otimes 3}\rar H^3(V, (\Omega^2_V)^{\otimes 3})$
with the map $\Tr\circ H^3(\gamma)$, where $\gamma$ is defined by
\begin{equation}\label{skew-cubic}
\gamma:(\Omega^2_V)^{\otimes 3}\rar \Omega^3_V\ ,\ \ \ 
\xi_1\otimes \xi_2\otimes \xi_3\mapsto
(s\alp^{-1}(\xi_1)\lrcorner\,\xi_2)\wedge \xi_3
\end{equation}
and $\Tr$ is the canonical isomorphism $H^3(V,\Omega^3_V)\isoto\CC$.
\end{proposition}

\begin{proof}
(i) The inclusion $s$ factors through $t:\TTT_V(-\log S)\into \TTT_V$,
the map identifying 
$\TTT_V(-\log S)$ as the subsheaf
of germs of vector fields tangent to $S$.
The fact that $H^1(t)$ coincides with the map, associating to
each first order deformation of the pair $(S,V)$ the respective first order
deformation of its second component $V$, is proved in the same way
as a similar property of $H^1(r)$ in \cite{B2}, Proposition 1.1. The wanted assertion then follows from the relation $s=t|_{\ker H^1(r)}$.

(ii) By part (i), Lemma \ref{bar-nabla} and formula (\ref{nabla}),
the Kodaira-Spencer class of $\pi$ at $b$ is the morphism
$e_b:T_bB=H^1(V,\Omega^2_V)\lra
\Hom(H^1(V,\Omega^2_V), H^2(V,\Omega^1_V))$ given by
$e_b:a\mapsto [v\mapsto \tilde a\lrcorner\, v]$, where
$\tilde a=H^1(s\circ \alp^{-1})(a)$. The Donagi--Markman cubic $Q$ is
the image of $e_b$ under the natural isomorphism
$\Hom (E,\Hom(E,E^\dual))\isoto (E\otimes E\otimes E)^\dual$,
where $E=H^1(V,\Omega^2_V)$ and $H^2(V,\Omega^1_V)$ is identified with
$E^\dual$ via the pairing $H^1(V,\Omega^2_V)\otimes
H^2(V,\Omega^1_V)\xrightarrow{\wedge} H^3(V,\Omega^3_V)=\CC$.
Hence $Q$ is the map
$E^{\otimes 3}\rar H^3(V,\Omega^3_V)=\CC\ , \ \ 
a\otimes v\otimes w\mapsto \mapsto (\tilde a\lrcorner\,v)\wedge w$,
which was to be proved.
\end{proof}

\begin{corollary}\label{quasi-symp}
Under the hypotheses and in the notation  of Proposition \ref{ebQ},
the Donagi--Markman cubic $Q$ is symmetric, so $\pi$ is a quasi-Lagrangian
fibration.
\end{corollary}

\begin{proof}
The fact that (\ref{skew-cubic}) is skew-symmetric in $\xi_i$
is an easy exercise in linear algebra. Then the induced pairing
on the cohomology $\otimes_{i=1}^3H^{p_i}(V,\Omega^2_V)\rar
H^{p_1+p_2+p_3}(V,\Omega^3_V)$ is graded skew-symmetric, that is,
the sign change resulting from the transposition of the first
and the second arguments is $(-1)^{(p_1+1)(p_2+1)}$, and similarly
for the transposition of the second and the third arguments.
It remains to apply Lemma \ref{weak}.
\end{proof}

\begin{remark}
If we replace $V$ by a Calabi--Yau 3-fold and choose for $\alp_V$
a generator of $H^0(V,\Omega^3_V)$, then the cubic form constructed
in part (ii) of Proposition \ref{ebQ} will be nothing else but
the Yukawa coupling, well known from the superstring theory
compactified on Calabi--Yau spaces.
\end{remark}

In order to get a genuine symplectic structure on $\JJJ$, we will have
to choose the isomorphism $j:\TTT_B^\dual\rar \EEE$ with more care. This choice
is naturally explained in terms of
the variation of the mixed Hodge structure on the cohomology
of $V\setminus S$.

\section{Variations of mixed Hodge structures}\label{VMHS}

A mixed Hodge structure (MHS) on a finite-dimensional $\CC$-vector space $H=H_\CC$
is the following set of data:
\begin{itemize}
\item[--] An integer lattice $H_\ZZ\subset H$ such that $H=H_\ZZ\otimes\CC$.
\item[--] A decreasing filtration $F^\fatdot H$ such that $F^0H=H$ and
$F^iH=0$ for $i\gg 0$ (Hodge filtration).
\item[--] An increasing filtration $W_\fatdot H$ defined over $\QQ$
(weight filtration)
such that $W_{-1}H=0$, $W_mH=H$ for $m\gg 0$, and such that $F^\fatdot$
induces on $\Gr_mH=W_mH/W_{m-1}H$ a pure Hodge structure of weight $m$
for any $m\geq 0$.
\end{itemize}

Here the induced Hodge filtration on $\Gr_mH$ is given by
$$F^p\Gr_mH=(F^pH\cap W_mH+W_{m-1}H)/W_{m-1}H.$$

Deligne \cite{De} introduced a natural MHS on the cohomology of a nonsingular
algebraic variety $Y$ depending functorially on $Y$. It is described as follows.
Let $\bar Y$ be a nonsingular complete variety containing $Y$ such that
$\bar Y\setminus Y=D$ is a divisor with simple normal crossings. Then
$H^\fatdot(Y,\CC)=\HHB^\fatdot(\bar Y,\Omega^\fatdot_{\bar Y}(\log D))$. The Hodge
and weight filtrations on the hypercohomology are defined via the respective filtrations on the complex of logarithmic differential forms. For the Hodge
filtration, we have
$$
F^p(\Omega^\fatdot_{\bar Y}(\log D))= \left(
\Omega^p_{\bar Y}(\log D)\xrightarrow{d}
\Omega^{p+1}_{\bar Y}(\log D)\xrightarrow{d}\ldots\ \  \right)[p],
$$
where $K^\fatdot[p]$ denotes the complex obtained from $K^\fatdot$ by the shift
$p$ steps to the right (so that $\Omega^p_{\bar Y}(\log D)$
is placed on the spot number $p$). 
The subcomplex
$
W_k(\Omega^\fatdot_{\bar Y}(\log D))\subset \Omega^\fatdot_{\bar Y}(\log D),
$
by definition, consists
of logarithmic forms with at most $k$ poles. Then
$$
F^pH^n(Y,\CC)=\im \HHB^n(F^p(\Omega^\fatdot_{\bar Y}(\log D)))\ ,\ \ 
W_{k+n}H^n(Y,\CC)=\im \HHB^n(W_k(\Omega^\fatdot_{\bar Y}(\log D))),
$$
the images being taken under the maps induced on the hypercohomology
by the natural inclusions of the complexes. It is proved that $W_\fatdot$
is defined over $\QQ$ and that the resulting MHS does not depend on the
completion $\bar Y$ of $Y$. Moreover, the spectral sequence of the
Hodge filtration on the complex $\Omega^\fatdot_{\bar Y}(\log D)$ converging to
the hypercohomology of the latter degenerates at $E_1$, so that
$$F^pH^n(Y,\CC)/F^{p+1}H^n(Y,\CC)\simeq H^{n-p}(\bar Y,
\Omega^p_{\bar Y}(\log D)).$$

In the particular case when $D$ is a smooth hypersurface in $Y$,
the weight filtration has only two non-trivial graded factors,
and the MHS is easily described via the Gysin exact sequence
\begin{equation}\label{gysin}
	\ldots\rar H^{q-2}(D)\xrightarrow{i_*} H^q(\bar Y)\xrightarrow{u^*} H^q(Y)
\xrightarrow{\Res}	 H^{q-1}(D)\rar H^{q+1}(\bar Y)\rar \ldots\ ,
\end{equation}
where $Y=\bar Y\setminus D$, the coefficients of the cohomology groups
are in $\ZZ$ or in a field, $i:D\into \bar Y$, $u:Y\into \bar Y$
are natural inclusions, $i_*$ is the Gysin homomorphism induced
by $i$, and $\Res$ is the residue (or the tube) map, see \cite{G}.
Here $W_{q-1}H^q(Y)=0$, $W_{q+1}H^q(Y)=H^q(Y)$,
$_{\scriptscriptstyle W}\!\gr_qH^q(Y)=W_qH^q(Y)=u^*H^q(\bar Y)$ is a pure Hodge structure of weight $q$,
a quotient of the one on $H^q(\bar Y)$,
and $_{\scriptscriptstyle W}\!\gr_{q+1}H^q(Y)=W_{q+1}H^q(Y)/W_qH^q(Y)\cong \im(\Res)$ is a pure Hodge substructure of $H^{q-1}(D)$ with weight shifted by 2.

A variation of mixed Hodge structures (VMHS) over a base $B$
is a quadruple $(\HHH,\HHH_\ZZ,\FFF^\fatdot,\WWW_\fatdot)$ such that:
\begin{itemize}
	\item[--] $\HHH$ is a holomorphic vector bundle on $B$.
	\item[--] $\HHH_\ZZ\subset \HHH$ is
a local system of free $\ZZ$-modules such that $\HHH_\ZZ\otimes\OOO_B
\simeq\HHH$.
  \item[--] $\FFF^\fatdot$ is a decreasing filtration
of $\HHH$ by holomorphic subbundles.
  \item[--] $\WWW_\fatdot$ is an increasing filtration of  $\HHH_\QQ=\HHH_\ZZ\otimes\QQ$ by local
systems of constant rank.
  \item[--] For each $b\in B$, the lattice $H_\ZZ:=\HHH_\ZZ|_b$ and the filtrations
$W_kH:=\WWW_k|_b$, $F^pH:=\FFF^p|_b$ define on $H:=\HHH|_b$ a MHS.
  \item[--] The Gauss--Manin connection $\nabla$ on $\HHH$ associated to the local
  system $\HHH_\ZZ\otimes\CC\subset\HHH$ satisfies the 
  transversality condition $\nabla (\FFF^p)\subset \FFF^{p-1}\otimes\Omega^1_B$.
\end{itemize}

Any family of pairs $(\bar Y,D)$ as above defines a VMHS
on the associated bundle of cohomology groups $H^n(\bar Y\setminus D,\CC)$.
By a family of pairs we mean a smooth proper morphism $f:\bar\YYY\rar B$
with connected fibers and a normal crossing divisor 
$\DDD=\DDD_1\cup\ldots\cup\DDD_r$ in $ \bar\YYY$
such that for any $k$-uple $(i_1,\ldots i_k)$ of distinct indices
from $\{1,\ldots, r\}$, the scheme-theoretic intersection $\DDD_{i_1}\cap\ldots\cap\DDD_{i_k}$,
if nonempty, is smooth over $B$ and is of codimension $k$ in $\bar \YYY$.
We denote the corresponding VMHS by $(\HHH^n(f,\DDD), \HHH^n_\ZZ(f,\DDD),
\FFF^\fatdot\HHH^n(f,\DDD),\WWW_\fatdot\HHH^n(f,\DDD))$.

As we have already mentioned in Section \ref{K3Ff}, the space of infinitesimal
deformations of a pair $(\bar Y,D)$ is canonically isomorphic to
$H^1(\bar Y, \TTT_{\bar Y}(-\log D))$, thus a family of pairs
defines the Kodaira--Spencer maps
$$
KS(f,\DDD)_b:T_bB\lra H^1(\bar Y_b, \TTT_{\bar Y_b}(-\log D_b))\ \ 
\mbox{for any}\ \ b\in B.
$$
If $h^1(\TTT_{\bar Y_b}(-\log D_b))$ is constant on $B$, then
the maps $KS(f,\DDD)_b$ fit into a morphism of vector bundles $KS(f,\DDD):\TTT_B\rar
R^1f_*(\TTT_{\bar\YYY/B}(-\log\DDD))$. Furthermore, there
is a natural contraction map
$$
\TTT_{\bar Y_b}(-\log D_b)\times
\Omega^p_{\bar Y_b}(\log D_b)\xrightarrow{\ \ \lrcorner\ \ } \Omega^{p-1}_{\bar Y_b}(\log D_b),
$$
and the $p$-th graded piece of the Gauss--Manin connection
on $\HHH^n(f,\DDD)$ with respect to the Hodge filtration 
$
\gr^p_F\nabla_b:(\FFF^p/\FFF^{p+1})|_b\rar(\FFF^{p-1}/\FFF^{p})\otimes\Omega^1_B|_b
$
is nothing but the contraction
$$
H^{n-p}(\bar Y_b,
\Omega^p_{\bar Y_b}(\log D_b))\xrightarrow{\ \ \lrcorner KS(f,\DDD)_b\ \ }
H^{n-p+1}(\bar Y_b,
\Omega^{p-1}_{\bar Y_b}(\log D_b))\otimes\Omega^1_B|_b
$$
with the class
$KS(f,\DDD)_b$, considered as an element of
$H^1(\bar Y_b, \TTT_{\bar Y_b}(-\log D_b)\otimes\Omega^1_B|_b)$.

\section{Symplectic structure}\label{SS}

We are resuming some of the notation of Sect. \ref{K3Ff}:
let $S$ be a K3 surface, $R\subset\Pic S$ an even integral lattice 
of signature $(1, \rk (R)-1)$ with a distinguished element
$\rho$, $\rho^2=2g-2$, and $B=\FFF^R_S$ the moduli stack of
Fano 3-folds $V$ with $(\Pic V, -K_V)\simeq (R,\rho)$ containing
$S$ as an anticanonical divisor, the
quadratic form on $\Pic V$ being defined by
$(D_1,D_2)\mapsto (-K_V\cdot D_1\cdot D_2)$.
Assume that $B$ is of dimension $>0$ and denote by
$f:\VVV\rar B$ the universal family over $B$.

\begin{theorem}
The relative Jacobian $\pi:\JJJ\rar B$ of the universal family
$f:\VVV\rar B$ is a Lagrangian fibration.
\end{theorem}

\begin{proof}
In order to compute the MHS on $H^3(Y)$ for $Y=V\setminus S$,
we apply \eqref{gysin} to the pair $(\bar Y,D)=(V,S)$:
\begin{equation}\label{ces}
	0\rar H^3(V)\xrightarrow{u^*} H^3(Y)
	\xrightarrow{\Res}	H^2(S)\xrightarrow{i_*} H^4(V)\rar 0.
\end{equation}
We will denote by $h^{p,q}$, $b_k$ the Hodge, resp. Betti
numbers of $V$.
The Hodge filtration on $H^3(Y)$ is given by
\begin{multline*}F^4=0,\ 
F^3=H^0(\Omega^3_V(S))\simeq\CC,\ 
F^2=F^3\oplus H^1(\Omega^2_V(\log S))\simeq \CC^{h^{2,1}+21-b_2},\\ 
F^1=F^2\oplus H^2(\Omega^1_V(\log S))=H^3(Y,\CC)\simeq \CC^{b_3+22-b_2}.
\end{multline*}
Further,
\begin{multline*}
W_3=H^3(V)\subset H^3(Y),\ W_3\cap F^3=0,\ W_3\cap F^2= H^1(\Omega^2_V),\\
W_3\cap F^1=H^1(\Omega^2_V)\oplus H^2(\Omega^1_V)=H^3(V),
\end{multline*}
and the Hodge structure of $_{\scriptscriptstyle W}\!\gr_4H^3(Y)$
coincides with that of $R^\perp=\ker i_*$ with weight shifted by 2, where the orthogonal complement is
taken in $H^2(S,\ZZ)$ with respect to the intersection form. Remark
that if $S$ is generic with the property $R\subset \Pic S$, 
then $R^\perp=(\Pic S)^\perp$ is nothing else but the lattice of transcendental cycles $\mathbf T_S$, so that $_{\scriptscriptstyle W}\!\gr_4H^3(Y)
\simeq \mathbf T_S(-1)$ as Hodge structures.

To get a symplectic structure on $\JJJ$, we will apply the criterion
of Lemma \ref{strong}.  We are using the same isomorphisms
$j_b:T_bB^\dual\rar T_0J(V_b)$ as in Proposition \ref{ebQ}, depending
on the choice of a generator $\alp_{V_b}$ of $H^0(\Omega^3_{V_b}(S))$.
To define
a global isomorphism $j:\TTT_B^\dual\rar \EEE$, where $\EEE$ is realized
as the restriction of $\TTT_{\JJJ/B}$ to the zero section of $\JJJ/B$,
we should fix some section $\alp_\VVV$ of the Hodge bundle $\FFF^3$.
We will pick up a generator $\alp_S$ of $H^0(\Omega^2_S)$ and
determine $\alp_\VVV$ by the requirement that $\Res(\alp_{V_b})=\alp_S$
for all $b\in B$; this provides a global trivialization of $\FFF^3$
over $B$.

For any local section $\gamma$ of $\HHH_3(\VVV/B,\ZZ)$, we want
$j^{-1}(\gamma)$ to be a closed 1-form on $B$, where $\HHH_3(\VVV/B,\ZZ)$
is naturally embedded into $\EEE=(F^2\HHH_3(\VVV/B))^\dual$ via the integration
map $\gamma\mapsto\int_\gamma$. The 1-form $j^{-1}(\gamma)$ is determined
by its values on the tangent vectors 
$a\in T_bB=H^1(V,\TTT_{V_b}(-S))$: 
$$j^{-1}(\gamma)(a)=\int_{\gamma_b}a\lrcorner\,\alp_{V_b}.
$$
If we try to mimic the proof of Theorem 7.7 in \cite{DM}, then in order
to verify the closedness of $j^{-1}(\gamma)$, we should find a function $g$, defined locally on $B$, whose
differential is $j^{-1}(\gamma)$, and the natural candidate
for $g$ would be $b\mapsto g(b)=\int_{\gamma_b}\alp_{V_b}$. The problem is
that $F^3H^3(Y)$ does not couple with $H_3(V)$ (here $V=V_b$, 
$Y=Y_b=V_b\setminus S$), so we have to correct the definition
of $g$ in lifting $\gamma_b$ to $H_3(Y)$.
Applying the Poincar\'e--Lefschetz duality to the cohomology exact
sequence of the pair $(V,S)$, we obtain
an exact sequence for the homology of $Y$ similar to \eqref{ces}:
$$
0\rar H_4(V)\rar H_2(S)\rar H_3(Y)\rar H_3(V)\rar 0.
$$
It shows that $\gamma$ can be lifted to a local section of $\HHH_3((\VVV
\setminus\mathbb S)/B,\ZZ)$, where $\mathbb S=S\times B\subset \VVV$. We
denote
the lifted section by the same symbol $\gamma$, then the above formula for $g$
provides a well-defined function on $B$. Its differential at a
point $b$ is computed in
terms of the Gauss--Manin connection of the VMHS associated to $f$:
$$
d_bg(a)=\int_{\gamma_b}\nabla_a\alp_\VVV.
$$
We have $\nabla_a\alp_\VVV\equiv a\lrcorner\,\alp_{V}\mod F^3H^3(Y)$ and
$a\lrcorner\,\alp_{V}\in H^1(\Omega^2_V)\subset W_3$. Since this
holds at any $b\in B$ and for any $a\in T_bB$, the subbundle
$E=\WWW_3\oplus\FFF^3\subset \HHH^3((\VVV\setminus\mathbb S)/B)$ is
$\nabla$-invariant. Picking up a flat basis of $\WWW_3$ and completing
it by $\alp_\VVV$ to a basis of $\WWW_3\oplus\FFF^3$, we obtain a basis of $E$
in which  the matrix of $\nabla|_E$ has the following form:
$$
C=\left[\begin{array}{ccc|c}
& & & c_{1n} \\
& \mbox{{\Huge 0}}& & \vdots \\
& & & c_{nn}
\end{array}\right]\ \ , \ \ \ c_{kn}\in H^0(\Omega^1_B).
$$
We have $j^{-1}(\gamma)=dg-gc_{nn}$. 
As $\Res (\alp_\VVV)$ is constant over $B$ and $\Res : H^3(V)\rar H^2(S)$
is defined over $\ZZ$, $c_{nn}=0$ and $j^{-1}(\gamma)$ is closed.

\end{proof}

\end{document}